\documentclass{commat}

\newcommand{\A}{{\mathcal A}}
\newcommand{\C}{{\mathcal C}}
\newcommand{\R}{{\mathcal R}}
\newcommand{\cee}{{\mathbb C}}
\newcommand{\que}{{\mathbb Q}}
\newcommand{\real}{{\mathbb R}}
\newcommand{\zed}{{\mathbb Z}}
\newcommand{\btau}{{\bar{\tau}}}
\newcommand{\be}{{\boldsymbol e}}
\newcommand{\bb}{{\boldsymbol b}}
\newcommand{\bv}{{\boldsymbol v}}
\newcommand{\bx}{{\boldsymbol x}}
\newcommand{\bwy}{{\boldsymbol y}}
\newcommand{\bo}{{\boldsymbol 0}}
\DeclareMathOperator{\spn}{span}
\DeclareMathOperator{\Tr}{Tr}

\title{%
    On average coherence of cyclotomic lattices
}

\author{%
    Lenny Fukshansky and David Kogan
}

\affiliation{
    \address{Lenny Fukshansky --
    Department of Mathematics, 850 Columbia Avenue, Claremont McKenna College, Claremont, CA 91711
}
    \email{%
    lenny@cmc.edu
}
    \address{David Kogan --
    Institute of Mathematical Sciences, Claremont Graduate University, Claremont, CA 91711
}
    \email{%
    david.kogan@cgu.edu
}
}

\abstract{%
    We introduce maximal and average coherence on lattices by analogy with these notions on frames in Euclidean spaces. Lattices with low coherence can be of interest in signal processing, whereas lattices with high orthogonality defect are of interest in sphere packing problems. As such, coherence and orthogonality defect are different measures of the extent to which a lattice fails to be orthogonal, and maximizing their quotient (normalized for the number of minimal vectors with respect to dimension) gives lattices with particularly good optimization properties. While orthogonality defect is a fairly classical and well-studied notion on various families of lattices, coherence is not. We investigate coherence properties of a nice family of algebraic lattices coming from rings of integers in cyclotomic number fields, proving a simple formula for their average coherence. We look at some examples of such lattices and compare their coherence properties to those of the standard root lattices.}



\VOLUME{31}
\YEAR{2023}
\NUMBER{1}
\firstpage{57}
\DOI{https://doi.org/10.46298/cm.10241}

\begin{paper}

\section{Introduction}\label{intro}

Let $L \subset \real^d$ be a lattice of full rank $d \geq 1$ in the Euclidean space $\real^d$, where we will always write $\|\ \|$ for the corresponding Euclidean norm. Define the (squared) \textit{minimum} of $L$ to be
\[
|L| : = \min \left\{\|\bx\|^2 : \bx \in L \setminus \{\bo \} \right\},
\]
and the set of \textit{minimal vectors} of $L$ to be
\[
S(L) : = \left\{\bx \in L : \|\bx\|^2 = |L| \right\}.
\]
The lattice $L$ is called \textit{well-rounded} (abbreviated WR) if $\spn_{\real} S(L) = \real^d$. There is a stronger condition for $L$ to be \textit{generated by minimal vectors} if $\spn_{\zed} S(L) = L$ (see~\cite{pohst}), and an even stronger condition for $L$ to have a \textit{basis of minimal vectors}, i.e. for $S(L)$ to contain a basis for $L$ (see~\cite{mart-schur}). We can associate a sphere packing to the lattice $L$ by placing maximal non-overlapping spheres of equal radius at the lattice points, then the radius of these spheres, called the \textit{packing radius} of $L$ will be $\sqrt{|L|}/2$ and the density of this lattice packing will be
\[
\delta(L) : = \frac{v_d |L|^{\frac{d}{2}}}{2^d \det(L)},
\]
where $v_d$ is the volume of a unit ball in~$\real^d$ and $\det(L)$ is the determinant of $L$, i.e., $\det(L) : = |\det(B)|$ for any choice of a $d \times d$ matrix $B$, called \textit{basis matrix} for $L$, such that $L = B \zed^d$. The determinant of $L$ is precisely the volume of any fundamental domain of $L$, such as the parallelepiped spanned by the column vectors of $B$. Given a basis matrix $B =
\begin{pmatrix}
\bb_1 & \dots & \bb_d
\end{pmatrix}
$ for $L$, we define the \textit{orthogonality defect} of $B$ as
\[
\nu(B) : = \frac{\prod_{j = 1}^d \|\bb_j\|}{\det(L)},
\]
i.e. the ratio of the volume of a rectangular box with sides $\|\bb_1\|,\dots,\|\bb_d\|$ to the volume of the parallelepiped spanned by the column vectors of $B$. Naturally, the Hadamard inequality $\nu(B) \geq 1$ holds with equality if and only if $B$ is an orthogonal basis. If $B \subseteq S(L)$, then
\begin{equation}\label{nu_delta}
\nu(B) = \frac{|L|^{\frac{d}{2}}}{\det(L)} = \frac{2^d}{v_d} \delta(L)
\end{equation}
is an invariant of the lattice $L$, which we will call the orthogonality defect of $L$ and denote by $\nu(L)$. Hence for a lattice with a basis of minimal vectors the packing density is proportionate to the orthogonality defect, i.e. to maximize the packing density one wants a lattice with a ``least orthogonal" minimal basis. Orthogonality defect figures prominently in lattice theory, especially in connection with algorithmic lattice problems (see~\cite{micciancio}). See also~\cite{conway:sloane} and~\cite{martinet} for detailed authoritative expositions of the theory of lattices and its fundamental connections to optimization problems, such as sphere packing and others.

Another measure of orthogonality for a collection of vectors is given by coherence and comes from signal processing. Given a finite set of vectors $S \subset \real^d$, we define its \textit{maximal coherence} as
\[
\C(S) : = \max \left\{\frac{\left| \left< \bx, \bwy \right> \right|}{\|\bx\| \|\bwy\|} : \bx \neq \bwy \in S \right\},
\]
where $\left<\,\ \right>$ stands for the usual Euclidean inner product, and its \textit{average coherence} as
\[
\A(S) : = \frac{1}{|S|-1} \max \left\{\sum_{\bwy \in S \setminus \{\bx\}} \frac{\left| \left< \bx, \bwy \right> \right|}{\|\bx\| \|\bwy\|} : \bx \in S \right\}.
\]
It is easy to see that $\A(S) = 0$ if and only if $S$ is an orthogonal collection of vectors, which in particular implies $|S| \leq d$. An important problem in signal processing is the construction of sufficiently large sets $S$ ($|S| > d$) with sufficiently low coherence. Special attention among such low-coherence sets is usually given to frames, which are over-determined spanning sets with certain additional properties, especially to the uniform tight frames: a finite set $S \subset \real^d$ is called a \textit{uniform tight frame} if all vectors in $S$ have the same norm and there exists a real constant $\gamma > 0$ such that
\[
\|\bv\| = \gamma \sum_{\bx \in S} \left< \bv, \bx \right>^2,
\]
for every $\bv \in \real^d$ (see~\cite{waldron} for a comprehensive exposition of tight frame theory).

We can extend the notion of coherence to lattices as follows. Notice that minimal vectors of a lattice $L$ come in $\pm$ pairs: $\bx \in S(L)$ if and only if $-\bx \in S(L)$. Then define $S'(L)$ to be a subset of $S(L)$ constructed by selecting one vector out of each such pair, and define maximal and average coherence of $L$ to be
\[
\C(L) : = \C(S'(L)),\ \A(L) : = \A(S'(L)),
\]
respectively. These values do not depend on the specific choice of vectors in $S'(L)$ out of each $\pm$ pair. If $L$ has a basis of minimal vectors, then $\A(L)$ becomes a certain alternative measure of its ``non-orthogonality": $\A(L) \geq 0$ with equality if and only if $S'(L)$ is an orthogonal basis for $L$. Maximal coherence on lattices has previously been introduced in~\cite{ejc-me} and studied on nearly orthogonal lattices in~\cite{david-me}, but average coherence has not previously been extended to lattices, as far as we know. Average coherence for frames was introduced in~\cite{bajwa}. Our definition of average coherence slightly differs from the one introduced in~\cite{bajwa}: in their definition, the absolute value is outside of the sum. We choose to move absolute value inside to ensure that the average coherence does not depend on the choice of the vectors in~$S'(L)$: it does not matter which vector from each $\pm$ pair in $S(L)$ is selected.

While there can be a relation between average coherence and orthogonality defect in some special cases, there does not appear to be a general dependence. On the other hand, it is interesting to understand which lattices with relatively large sets of minimal vectors simultaneously have small average coherence and large orthogonality defect. To this end, given a lattice $L \subset \real^d$ with a basis of minimal vectors, we define its \textit{orthogonality product measure} (referred to from here on simply as \textit{product measure}) to be
\begin{equation}\label{prod_measure}
\Pi(L) : = \frac{|S'(L)| \nu(L)}{d \A(L)}.
\end{equation}
Then a lattice $L$ with large $|S'(L)|$ (as compared to the dimension $d$), small $\A(L)$ and large $\nu(L)$ will have large $\Pi(L)$. We can then ask which lattices have large $\Pi(L)$. In this note, we investigate average coherence and product measure on the family of cyclotomic lattices, a special family of ideal lattices. We start out by introducing the ideal lattices.

Let $K$ be a number field of degree $d$ over $\que$, and let $\mathcal O_K$ be its ring of integers. Let
\[
\sigma_1,\dots,\sigma_{r_1},\tau_1,\bar{\tau}_1,\dots,\tau_{r_2},\bar{\tau}_{r_2} : K \hookrightarrow \cee
\]
be its embeddings into the field of complex numbers, where $r_1+2r_2 = d$ and $\sigma_1,\dots,\sigma_{r_1}$ are real embeddings, whereas $\tau_1,\bar{\tau}_1,\dots,\tau_{r_2},\bar{\tau}_{r_2}$ are pairs of complex conjugate embeddings. The \textit{Minkowski embedding} of $K$ into $\real^d$ is then defined as
\[
\Sigma_K : = \left(\sigma_1,\dots,\sigma_{r_1},\Re(\tau_1),\Im(\tau_1),\dots,\Re(\tau_{r_2}),\Im(\tau_{r_2})\right) : K \hookrightarrow \real^d,
\]
and the image of $\mathcal O_K$ under this embedding, $\Lambda_K : = \Sigma_K(\mathcal O_K)$ is a Euclidean lattice of full rank in~$\real^d$. Furthermore,
\begin{equation}\label{det_L}
\det(\Lambda_K) = 2^{-r_2} |\Delta_K|^{1/2},
\end{equation}
where $\Delta_K$ stands for the discriminant of $K$. Such lattices are called \textit{number field lattices}; they form a special case of the more general \textit{ideal lattices} (of trace type), which are given by the same construction on an arbitrary fractional ideal in $K$. This construction of ideal lattices is classical: it can be found, for instance, in~\cite{bor:sha} (pp. 94--99) or~\cite{tsfasman} (Chapter 5.3), as well as in~\cite{ideal_lattices}.

We focus specifically on cyclotomic fields. Let $\zeta_n = e^{\frac{2 \pi i}{n}}$ for $n > 2$ be the $n$-th primitive root of unity and $K = \que(\zeta_n)$ be the corresponding $n$-th cyclotomic number field, then $d = [K:\que] = \phi(n)$ and the ring of integers $\mathcal O_K = \zed[\zeta_n]$. Then the group of $n$-th roots of unity
\[
\R_n : = \left\{\zeta_n^k : 1 \leq k \leq n \right\}
\]
is precisely the set of all roots of unity contained in $\mathcal O_K$. We refer to the lattice $\Lambda_K$ as the \textit{$n$-th cyclotomic lattice}. We give a more detailed description of cyclotomic lattices and their properties in Section~\ref{cyclotomic}, in particular explaining that they have bases of minimal vectors and
\[
|S'(\Lambda_K)| = \left\{
\begin{array}
{ll}
n & \mbox{if $n$ is odd,} \\
\frac{1}{2} n & \mbox{if $n$ is even.}
\end{array}
\right.
\]
Further, we demonstrate the well-known fact that in the cyclotomic case the orthogonality defect
\begin{equation}\label{nu_L}
\nu(\Lambda_K) = {\left(\frac{\phi(n)}{\prod_{p \mid n} p^{e_p - \frac{1}{p-1}}}\right)}^{\frac{\phi(n)}{2}},
\end{equation}
where $n = \prod_{p \mid n} p^{e_p}$ and the product in the denominator is over all primes dividing~$n$. We also define the average coherence $\A(\alpha)$ for any $\alpha \in S'(\Lambda_K)$, as well as $\A(\Lambda_K)$, the average coherence of the lattice $\Lambda_K$, in~\eqref{A_alpha} and~\eqref{A_LK}, respectively. Finally, cyclotomic lattices are \textit{strongly eutactic}, meaning that their sets of minimal vectors form uniform tight frames in their respective Euclidean spaces.

Cyclotomic lattices have been extensively studied in the context of lattice theory (see Section~8.7 of~\cite{conway:sloane} and references therein), and their structure is generally understood. One goal of this note is to attract some attention to the notions of average and maximal coherence on lattices. We use cyclotomic lattices as a simple and attractive case study. As it turns out, there is a particularly simple and elegant arithmetic formula for the average coherence of this family of lattices.

\begin{theorem}\label{main_cyclo} Let $n > 2$ be an integer, and let $\Lambda_K$ be the corresponding cyclotomic lattice for $K = \que(\zeta_n)$. Then
\[
\C(\Lambda_K) = \left\{
\begin{array}
{ll}
0 & \mbox{if $n$ is a power of $2$,} \\
\frac{1}{p-1} & \mbox{if $p$ is the smallest odd prime dividing $n$.}
\end{array}
\right.
\]
Additionally, for any $\alpha \in S'(\Lambda_K)$,
\begin{equation}\label{av_coh_formula}
\A(\alpha) = \A(\Lambda_K) = \left\{
\begin{array}
{ll}
\frac{2^{\omega(n)}-1}{n-1} & \mbox{if $n$ is odd,} \\
\frac{2^{\omega(n)}-2}{n-2} & \mbox{if $n$ is even,}
\end{array}
\right.
\end{equation}
where $\omega$ is the number of prime divisors function. Combining~\eqref{av_coh_formula} with \eqref{nu_L}, we readily obtain an explicit formula for $\Pi(\Lambda_{\que(\zeta_n)})$, which depends only on~$n$:
\[
\Pi(\Lambda_{\que(\zeta_n)}) = \left\{
\begin{array}
{ll}
\frac{n (n-2) \phi{(n)}^{\frac{\phi(n)}{2} - 1}}{2 (2^{\omega(n)} - 2) {\left(\prod_{p \mid n} p^{e_p - \frac{1}{p-1}}\right)}^{\frac{\phi(n)}{2}}} & \mbox{if $2 \mid n$}, \\
\frac{n (n-1) \phi{(n)}^{\frac{\phi(n)}{2} - 1}}{(2^{\omega(n)} - 1) {\left(\prod_{p \mid n} p^{e_p - \frac{1}{p-1}}\right)}^{\frac{\phi(n)}{2}}} & \mbox{if $2 \nmid n$}.
\end{array}
\right.
\]
\end{theorem}

We prove Theorem~\ref{main_cyclo} in Section~\ref{coherence}. In Section~\ref{coh_orth}, we demonstrate several examples, aiming to determine values of~$n$ for which $\Pi(\Lambda_{\que(\zeta_n)})$ is the largest in a fixed dimension $d = \phi(n)$. For comparison purposes, we also compute the coherence and product measure values for the standard root lattices. Of course, it is easy to see that the product measure values for cyclotomic lattices are not nearly as large as for the root lattices in the same dimensions. On the other hand, root lattices are truly exceptional (in particular, they are local maxima of the packing density function in their dimensions; see, for instance, Chapter~4 of~\cite{martinet} for details), and there are very few of them. Cyclotomic lattices present a larger family of lattices with interesting properties (in even dimensions given by the values of the Euler $\phi$-function), including numerous examples of lattices with low maximal coherence. In fact, as we discuss at the end of Section~\ref{coh_orth}, the maximal and average coherence of cyclotomic lattices, in contrast with the root lattices, are about the same on the average as $n \to \infty$, which can also make them potentially interesting from the standpoint of sparse signal processing: it guarantees that signal frequencies represented by the minimal vectors are well spread out, a useful feature for signal recovery (see~\cite{bajwa0} and~\cite{bajwa}). For future research, it would be interesting to investigate average coherence of other families of lattices coming from algebraic constructions, including some more general ideal lattices, as well as to study properties and general behavior of average coherence as a function on lattices. We are now ready to proceed.

\section{Cyclotomic lattices}\label{cyclotomic}

In this section we give an alternative and for our purposes more convenient description of cyclotomic lattices. Let $K = \que(\zeta_n)$ for $n > 2$, then $K$ only has the complex embeddings
\[
\tau_1, \btau_1, \dots, \tau_{d/2}, \btau_{d/2} : K \hookrightarrow \cee,
\]
so $r_1 = 0$ and $d = \phi(n) = 2r_2$. For each $\alpha \in \mathcal O_K$, the trace of $\alpha$ is given by
\[
\Tr_K(\alpha) : = \sum_{k = 1}^{d/2} (\tau_k(\alpha) + \btau_k(\alpha)).
\]
Using the notation of Section~8.7 of~\cite{conway:sloane} (see also~\cite{eva_bayer}), we can think of the cyclotomic lattice $\Lambda_K$ as the free $\zed$-module $\mathcal O_K$ equipped with the bilinear form
\[
\left< \alpha, \beta \right> : = \frac{1}{2} \Tr_K(\alpha \bar{\beta})
\]
for any $\alpha, \beta \in \mathcal O_K$. It is easy to verify that $\left< \alpha, \beta \right>$ is equal to the usual dot product of the vectors $\Sigma_K(\alpha)$ and $\Sigma_K(\beta)$ in $\real^d$. Then for any $\alpha = a + bi \in \mathcal O_K = \zed[\zeta_n]$ we have $\alpha \bar{\alpha} = a^2 + b^2$, and hence
\[
\left< \alpha, \alpha \right> = \sum_{k = 1}^{d/2} \tau_k(a^2 + b^2) = \sum_{k = 1}^{d/2} \left(\Re{(\tau_k(\alpha))}^2 + \Im{(\tau_k(\alpha))}^2 \right).
\]
By the results of~\cite{kate-me}, $\Lambda_K$ is WR with $|\Lambda_K| = \frac{\phi(n)}{2}$ and $\alpha \in S(\Lambda_K)$ if and only if it is a root of unity, i.e.
\[
S(\Lambda_K) = \{\pm \alpha : \alpha \in \R_n \} = \left\{
\begin{array}
{ll}
\R_n & \mbox{if $2 \mid n$} \\
\R_{2n} & \mbox{if $2 \nmid n$},
\end{array}
\right.
\]
since
\[
-1 = e^{\pi i} = \left\{
\begin{array}
{ll}
e^{\frac{2 (n/2) \pi i}{n}} & \mbox{if $2 \mid n$} \\
e^{\frac{2 n \pi i}{2n}} & \mbox{if $2 \nmid n$}.
\end{array}
\right.
\]
Let $\alpha, \beta \in S(\Lambda_K)$, then
\[
\left< \alpha, \beta \right> = \frac{1}{2} \Tr_K(\alpha \bar{\beta}),
\]
where $\alpha \bar{\beta}$ is also a root of unity. Suppose that $\alpha \bar{\beta}$ is an $m$-th primitive root of unity of for some $m \mid n$; then it is a root of $m$-th cyclotomic polynomial $\Phi_m(x)$. Notice that the trace of an algebraic number is the negative of the second coefficient of its minimal polynomial. It is a well-known fact that
\[
\Phi_m(x) = x^{\phi(m)} - \mu(m) x^{\phi(m)-1} + \dots,
\]
where $\mu$ is the M\"obius function. Hence $\Tr_{\que(\alpha \bar{\beta})}(\alpha \bar{\beta}) = \mu(m)$, and therefore
\begin{equation}\label{inn_ab}
\left< \alpha, \beta \right> = \frac{1}{2} \Tr_K(\alpha \bar{\beta}) = \frac{[K:\que(\alpha \bar{\beta})]}{2} \Tr_{\que(\alpha \bar{\beta})}(\alpha \bar{\beta}) = \frac{\phi(n)}{2 \phi(m)} \mu(m).
\end{equation}
Further, if $\alpha = \zeta_n^{k_1}$ and $\beta = \zeta_n^{k_2}$, then $m = \frac{n}{\gcd(k_1 - k_2, n)}$, and so the cosine of the angle between these two vectors is
\begin{equation}\label{c_ab}
c(\alpha,\beta) : = \frac{\left< \alpha, \beta \right>}{\sqrt{\left< \alpha, \alpha \right> \left< \beta, \beta \right>}} = \frac{\phi(n)}{\phi(n) \phi(m)} \mu(m) = \frac{\mu \left(\frac{n}{\gcd(k_1 - k_2, n)}\right)}{\phi \left(\frac{n}{\gcd(k_1 - k_2, n)}\right)}.
\end{equation}
Define $s : = |S(\Lambda_K)|$, so $s = n$ if $n$ is even and $s = 2n$ if $n$ is odd. Then we can write
\[
S(\Lambda_K) = \left\{\zeta_n^k, \zeta_n^{k+ \frac{s}{2}} : 1 \leq k \leq s/2 \right\},
\]
where $\zeta_n^{k+ \frac{s}{2}} = -\zeta_n^k$ and $c \left(\zeta_n^k, \zeta_n^{k+ \frac{s}{2}}\right) = -1$, as expected.
Hence let
\[
S'(\Lambda_K) = \left\{\zeta_n^k : 1 \leq k \leq s/2 \right\},
\]
so the coherence of the lattice $\Lambda_K$ is given by
\[
\C(\Lambda_K) = \max \left\{\left| c(\alpha,\beta) \right| : \alpha, \beta \in S'(\Lambda_K), \alpha \neq \beta \right\}.
\]
Then for any two $\alpha = \zeta^{k_1}, \beta = \zeta^{k_2} \in S'(\Lambda_K)$, $|k_1-k_2| \leq s/2 - 1$, so $c(\alpha,\beta) \neq \pm 1$. Additionally, for each $\alpha \in S'(\Lambda_K)$ define its average coherence to be
\begin{equation}\label{A_alpha}
\A(\alpha) = \frac{1}{|S'(\Lambda_K)| - 1} \sum_{\beta \in S'(\Lambda_K) \setminus \{\alpha \}} |c(\alpha,\beta)|.
\end{equation}
The average coherence of $\Lambda_K$ is then given by
\begin{equation}\label{A_LK}
\A(\Lambda_K) = \max \{\A(\alpha) : \alpha \in S'(\Lambda_K) \}.
\end{equation}
Now, the discriminant of the cyclotomic field $K = \que(\zeta_n)$ is given by
\[
\Delta_K = {(-1)}^{\frac{\phi(n)}{2}} n^{\phi(n)} \prod_{p \mid n} p^{- \frac{\phi(n)}{p-1}},
\]
where the product is over all primes $p$ dividing $n$ (see, for instance, Section~8.7.3 of~\cite{conway:sloane}). Combining this observation with~\eqref{nu_delta},~\eqref{det_L}, and the fact that $|\Lambda_K| = \phi(n)/2$, we obtain~\eqref{nu_L}.

We also briefly comment on the structure of cyclotomic lattices, which is well known (see, for instance, Section~8.7 of~\cite{conway:sloane}). Two lattices $L_1, L_2 \subset \real^k$ are called \textit{similar}, denoted $L_1 \sim L_2$, if there exists a nonzero real constant $\gamma$ and a $k \times k$ real orthogonal matrix $U$ such that $L_2 = \gamma U L_1$; if $\gamma = \pm 1$, $L_1$ and $L_2$ are \textit{isometric}, denoted $L_1 \cong L_2$. For any lattice $L \subset \real^d$ of rank $d$, its \textit{dual} is the lattice
\[
L^* : = \left\{\bx \in \real^d : \left< \bx, \bwy \right> \in \zed\ \forall\ \bwy \in L \right\}.
\]
The root lattice $A_n$ is defined as
\begin{equation}\label{A_n}
A_n = \left\{\bx \in \zed^{n+1} : \sum_{i = 1}^{n+1} x_i = 0 \right\},
\end{equation}
which is a lattice of rank $n$, as is its dual $A_n^*$. With this notation, the following is true:

\begin{enumerate}

\item If $n = p$ is an odd prime, then $\Lambda_K \sim A^*_{p-1}$,

\item If $n = p^k$ is an odd prime power, then $\Lambda_K \sim \bigoplus_{j = 1}^{p^{k-1}} A^*_{p-1}$,

\item If $n = p^k q^l$ is a product of two distinct odd prime powers, then
\[
\Lambda_K \sim \left(\bigoplus_{j = 1}^{p^{k-1}} A^*_{p-1}\right) \otimes \left(\bigoplus_{j = 1}^{q^{l-1}} A^*_{q-1}\right).
\]

\end{enumerate}

Lattices $A^*_n$ are known to be strongly eutactic. Further, tensor products of strongly eutactic lattices as well as direct sums of isometric strongly eutactic lattices are strongly eutactic (see Chapter~3 of~\cite{martinet}). This observation, along with the above properties, implies that cyclotomic lattices are in general strongly eutactic.

\section{Coherence of cyclotomic lattices}\label{coherence}

In this section we prove Theorem~\ref{main_cyclo} in a series of several lemmas. Throughout this section, $K = \que(\zeta_n)$ for the specified choices of $n$ and $\Lambda_K$ is the corresponding cyclotomic lattice.

\begin{lemma}\label{power-2} Suppose $n = 2^m$ for some $m \geq 1$, then $\Lambda_K$ is an orthogonal lattice, which is similar to~$\zed^{2^{m-1}}$. In particular, $\C(\Lambda_K) = 0$.
\end{lemma}

\proof
First notice that $\phi(2^m) = 2^{m-1}$, thus $\Lambda_K$ is a lattice of rank~$2^{m-1}$ with $2^m$ minimal vectors. Let $\alpha, \beta \in S'(\Lambda_K)$ and suppose $\alpha \bar{\beta}$ is a $k$-th primitive root of unity for some $k \mid 2^m$. Then $k = 2^l$ for some $0 \leq l \leq m$, and by~\eqref{inn_ab},
\[
\left< \alpha, \beta \right> = \frac{1}{2} \frac{\phi(2^m)}{\phi(2^l)} \mu(2^l) = 0,
\]
unless $l = 0$ or $1$. If $l = 0,1$, then $\alpha \bar{\beta}$ is either a first or second root of unity, i.e. $\alpha \bar{\beta} = \pm 1$, which implies that $\alpha = \pm \beta$. Therefore $c(\alpha,\beta) = 0$ for any pair of distinct minimal vectors in~$S'(\Lambda_K)$, and so~$S(\Lambda_K)$ consists of $\phi(n) = n/2 = 2^{m-1}$ plus-minus pairs of orthogonal basis vectors of equal norm. Hence $\Lambda_K \sim \zed^{2^{m-1}}$.
\endproof

\begin{lemma}\label{coh} Assume that $n$ is not a power of $2$, and let $p$ be the smallest odd prime dividing $n$. Then
\[
\C(\Lambda_K) = \frac{1}{p-1}.
\]
\end{lemma}

\proof
Let $\alpha = \zeta_n^{k_1} \in S'(\Lambda_K)$, then
\[
\left\{\beta \in S'(\Lambda_K) : \beta \neq \alpha \right\} = \left\{\zeta_n^{k_2} : 1 \leq k_2 \leq s/2, k_2 \neq k_1 \right\},
\]
and so $k_1-k_2$ takes on all nonzero integer values between $k_1-1$ and $k_1 - s/2$. In particular, $k_1-k_2 < s/2$, which means that $c(\alpha,\beta) \neq \pm 1$. Since $p$ is the smallest odd prime dividing $n$, $2 < p \leq s/2$. Then let $k_1 = p+1$ and $k_2 = 1$, and for the corresponding $\alpha = \zeta_n^{k_1}$, $\beta = \zeta_n^{k_2}$,~\eqref{c_ab} gives
\[
| c(\alpha,\beta) | = \frac{1}{p-1}.
\]
On the other hand, $\frac{n}{\gcd(k_1 - k_2, n)} \neq 1$ is a divisor of $n$, which cannot be equal to~$2$: in fact, notice that, if $\frac{n}{\gcd(k_1 - k_2, n)} = 2$, then $n$ is even and $|k_1-k_2| = n/2 = s/2$, however we know that $|k_1-k_2| \leq s/2 - 1$. Hence it cannot be smaller than $p$, and so~\eqref{c_ab} guarantees that $\C(\Lambda_K) \leq \frac{1}{p-1}$. Thus we have the result.
\endproof

\begin{lemma}\label{av_coh-odd} Assume $n$ is odd and square-free, then
\[
\A(\Lambda_K) = \frac{\tau(n)-1}{n-1},
\]
where $\tau(n)$ is the number of divisors of $n$.
\end{lemma}

\proof
Since $n$ is odd, we have $s/2 = n$. Let $\alpha = \zeta_n^k \in S'(\Lambda_K)$ for some $1 \leq k \leq s/2$, then by~\eqref{c_ab},
\begin{align*}
    \A(\alpha) 
    { }&{ }= \frac{1}{s/2-1} \sum_{j = 1,\ j \neq k}^{s/2} \frac{1}{\phi \left(\frac{n}{\gcd(j-k, n)}\right)} \\
    &{ }= \frac{1}{n-1} \sum_{m = 1-k,\ m \neq 0}^{n-k} \frac{1}{\phi \left(\frac{n}{\gcd(m, n)}\right)} \\
    &{ }= \frac{1}{n-1} \sum_{d \mid n,\ d \neq n} \frac{a_d}{\phi(n/d)},
\end{align*}
where $a_d = $ the number of times $\gcd(m,n) = d$ for nonzero $1-k \leq m \leq n-k$. Notice that the set $\{1-k,\dots,n-k \}$ is a complete residue system modulo~$n$, as is the set $\{0,\dots,n \}$ and hence the number of times $\gcd(m,n) = d$ for nonzero $1-k \leq m \leq n-k$ equals the number of times $\gcd(m,n) = d$ for $1 \leq m \leq n$. Therefore we can write
\[
\A(\alpha) = \frac{1}{n-1} \sum_{d \mid n, d \neq n} \frac{a_d}{\phi(n/d)},
\]
where
\[
a_d = \left| \left\{1 \leq m \leq n : \gcd(m,n) = d \right\} \right| = \phi(n/d),
\]
which is independent of $k$ and thus of the choice of $\alpha$. Hence we have
\[
\A(\Lambda_K) = \frac{1}{n-1} \sum_{d \mid n,\ d \neq n} \frac{\phi(n/d)}{\phi(n/d)} = \frac{1}{n-1} \sum_{d \mid n,\ d \neq n} 1 = \frac{\tau(n)-1}{n-1}.
\qedhere
\]
\endproof

\begin{lemma}\label{av_coh-even} Assume $n$ is even and square-free, then
\[
\A(\Lambda_K) = \frac{\tau(n)-2}{n-2},
\]
where $\tau(n)$ is the number of divisors of $n$.
\end{lemma}

\proof
Since $n$ is even, we have $s/2 = n/2$. Let $\alpha = \zeta_n^k \in S'(\Lambda_K)$ for some $1 \leq k \leq s/2$, then by~\eqref{c_ab},
\begin{align*}
    \A(\alpha)
    { }&{ }= \frac{1}{s/2-1} \sum_{j = 1,\ j \neq k}^{s/2} \frac{1}{\phi \left(\frac{n}{\gcd(j-k, n)}\right)} \\
    &{ }= \frac{2}{n-2} \sum_{m = 1-k,\ m \neq 0}^{\frac{n}{2} - k} \frac{1}{\phi \left(\frac{n}{\gcd(m, n)}\right)} \\
    &{ }= \frac{2}{n-2} \sum_{d \mid n,\ d < \frac{n}{2}} \frac{b_d}{\phi(n/d)},
\end{align*}
where $b_d = $ the number of times $\gcd(m,n) = d$ for nonzero $1-k \leq m \leq \frac{n}{2} - k$. Notice that, if $d \neq 1,2$, then for any such $m$ there is a unique $m' = m+n/2$ such that
\[
    \gcd(m',n) = \gcd(m,n) = d
    \qquad \textup{and} \qquad
    \frac{n}{2} - k \leq m' \leq n-k.
\]
Therefore for each divisor $d \neq 1,2$ of $n$ with $d < n/2$, $b_d = \frac{\phi(n/d)}{2}$. On the other hand,
\[
\gcd(m,n) = 1 \Leftrightarrow \gcd(m',n) = 2,\ \gcd(m,n) = 2 \Leftrightarrow \gcd(m',n) = 1,
\]
so $b_1 + b_2 = \phi(n) = \phi(n/2)$. Further, observe that $d \mid n$ with $d < n/2$ if and only if $d \mid \frac{n}{2}$ and $d \neq n/2$. Hence
\begin{align*}
    \A(\Lambda_K)
    { }&{ }= \frac{2}{n-2} \left(\frac{\phi(n/2)}{\phi(n/2)} + \sum_{d \mid n,\ d < \frac{n}{2},\ d \neq 1,2} \frac{\phi(n/d)}{2 \phi(n/d)}\right) \\
    &{ }= \frac{2}{n-2} \left(1 + \frac{1}{2} (\tau(n)-4)\right) \\
    &{ }= \frac{\tau(n)-2}{n-2},
\end{align*}
since the number of divisors $d$ of $n$ such that $d < n/2$ is $\tau(n) -2$: we count all the divisors except for $n$ and $n/2$.
\endproof

\begin{corollary}\label{sqr} Let $n > 2$ be an integer and let $n' = \prod_{p \mid n} p$ be its square-free part. Let $\Lambda_K$ be the corresponding cyclotomic lattice for $K = \que(\zeta_n)$. Then for any $\alpha \in S'(\Lambda_K)$,
\[
\A(\alpha) = \A(\Lambda_K) = \left\{
\begin{array}
{ll}
\frac{\tau(n')-1}{n-1} & \mbox{if $n$ is odd,} \\
\frac{\tau(n')-2}{n-2} & \mbox{if $n$ is even.}
\end{array}
\right.
\]
\end{corollary}

\proof
For each $\alpha = \zeta_n^k \in S'(\Lambda_K)$, we have
\[
\A(\alpha) = \frac{1}{|S'(\Lambda_K)| - 1} \sum_{\beta \in S'(\Lambda_K) \setminus \{\alpha \}} |c(\alpha,\beta)| = \frac{2}{s-2} \sum_{j = 1,\ j \neq k}^{s/2} \frac{\left| \mu \left(\frac{n}{\gcd(j-k, n)}\right) \right|}{\phi \left(\frac{n}{\gcd(j-k, n)}\right)},
\]
where for each $\beta = \zeta_n^j \in S'(\Lambda_K)$,
\[
c(\alpha,\beta) = \frac{\mu \left(\frac{n}{\gcd(k - j, n)}\right)}{\phi \left(\frac{n}{\gcd(k-j, n)}\right)} = 0,
\]
unless $\frac{n}{\gcd(k-j, n)}$ is square-free, i.e. a divisor of~$n'$. Thus
\begin{equation}\label{cor-coh-1}
\A(\alpha) = \frac{2}{s-2} \sum_{m = 1-k,\ m \neq 0}^{\frac{s}{2}-k} \frac{\left| \mu \left(\frac{n}{\gcd(m, n)}\right) \right|}{\phi \left(\frac{n}{\gcd(m, n)}\right)} = \frac{2}{s-2} \sum_{\frac{n}{d} \mid n',\ d < \frac{s}{2}} \frac{c_d}{\phi(n/d)},
\end{equation}
where
\[
c_d = \left| \left\{1-k \leq m \leq \frac{s}{2}-k : \gcd(m,n) = d \right\} \right|.
\]
Notice that every divisor $d$ of $n$ such that $n/d$ divides $n'$ is of the form $d = d' (n/n')$, where $d' \mid n'$. Let $s' = n'$ if $n'$ is even and $2n'$ if $n'$ is odd, then
\[
c_d = \left| \left\{1-k \leq m \leq \frac{s'}{2}-k : \gcd(m,n') = d' \right\} \right| = \left\{
\begin{array}
{ll}
a_{d'} & \mbox{if $2 \nmid n'$} \\
b_{d'} & \mbox{if $2 \mid n'$},
\end{array}
\right.
\]
where $a_{d'}$ and $b_{d'}$ are as in Lemmas~\ref{av_coh-odd} and~\ref{av_coh-even}, respectively. The result then follows by combining~\eqref{cor-coh-1} with these lemmas.
\endproof

\proof[Proof of Theorem~\ref{main_cyclo}]
Notice that for any positive integer $n$ with its square-free part $n'$, $\tau(n') = 2^{\omega(n)}$. The statement of the theorem now follows by combining Lemmas~\ref{power-2},~\ref{coh} with Corollary~\ref{sqr}.
\endproof

\section{Coherence and orthogonality defect}\label{coh_orth}

Throughout this section, let us write $\C_n$, $\A_n$, $\nu_n$ and $\Pi_n$ for $\C(\Lambda_{\que(\zeta_n)})$, $\A(\Lambda_{\que(\zeta_n)})$, $\nu(\Lambda_{\que(\zeta_n)})$, and $\Pi(\Lambda_{\que(\zeta_n)})$, respectively. We aim to understand the behavior of these functions as $n$ ranges through natural numbers. The first observation is that for odd $n$, $\Lambda_{\que(\zeta_{2n})} = \Lambda_{\que(\zeta_n)}$, and the formulas from Section~\ref{intro} yield
\[
\C_{2n} = \C_n,\ \A_{2n} = \A_n,\ \nu_{2n} = \nu_n,\ \Pi_{2n} = \Pi_n,
\]
as expected.

Let us start by briefly recalling the order of the arithmetic function $\phi(n)$ (see Chapter~18 of~\cite{hardy_wright} for further details). For all $n > 2$,
\begin{equation}\label{phi_ineq-1}
\frac{n}{e^{\gamma} \log \log n + \frac{3}{\log \log n}} < \phi(n) < n,
\end{equation}
where $\gamma = 0.57721\dots$ is Euler's constant. In fact, $\phi(n) < \frac{n}{e^{\gamma} \log \log n}$ for infinitely many $n$, although the average order of $\phi(n)$ is
\[
\frac{1}{n} \sum_{m = 1}^n \phi(m) = \frac{3n}{\pi^2} + O(\log n).
\]
Recall now that $s/2$, the cardinality of $S'(\Lambda_{\que(\zeta_n)})$ is $n$ or $n/2$, depending on the parity of $n$, whereas the rank of $\Lambda_{\que(\zeta_n)}$ is $\phi(n)$. Since it is desirable to have the number of minimal vectors as large as possible, compared to the dimension, we may want to consider values of $n$ for which $\phi(n)$ is close to the lower bound of~\eqref{phi_ineq-1}.

A particularly interesting situation from the stand-point of signal processing and of lattice theory arises when $2\phi(n)/s$ and $\A_n$ are small, while $\nu_n$ is large: this would mean that $S(\Lambda_{\que(\zeta_n)})$ is a configuration of many (in comparison to dimension) vectors, which are incoherent and non-orthogonal. Such configurations can be useful, for instance, in recovering signals transmitted with erasures (see~\cite{paulsen}). To this end, we observe that the values of $n$ that maximize $\Pi_n$ for each fixed dimension~$\phi(n)$ are large $n$ with small prime factors and small prime factor powers, and similarly for maximizing~$\nu_n$. On the other hand, values of $n$ minimizing $\A_n$ are large $n$ (for a fixed value of $\phi(n)$) with few prime factors, whereas $\C_n$ is minimized by $n$ with large prime factors. In particular, it appears that large $\Pi_n$ is more correlated with large $\nu_n$ than with small~$\A_n$. Indeed, consider the examples in Table~\ref{coh_orth_ex}: the values marked in bold are maximal among all $n$ with that value of $\phi(n)$ for $\nu_n$ and $\Pi_n$, and minimal for $\C_n$ and $\A_n$. We have also computed many additional examples, and the same observations seem to hold.

Further, although there is a general positive correlation between $\A_n$ and $\nu_n$ (see for instance dimension $24$ in Table~\ref{coh_orth_ex}), there are nevertheless sequences of closely related values of $n$ where the correlation is negative. Observe, for instance, dimension $72$ in Table~\ref{coh_orth_ex}. Take $n \in \left\{111, 117, 135, 228, 252\right\}$. If we arrange these in order of number of minimal vectors of $\Lambda_{\que(\zeta_n)}$, we have $s \in \left\{222, 228, 234, 252, 270\right\}$. These lattices respectively have $\A_n$ values of $0.0\overline{27}, 0.0265\dots, 0.0259\dots, 0.024$, and $0.0224\dots$. However, as $\A_n$ decreases, we see an increase in $\nu_n$, from $\nu_{111} = \nu_{222} = 2447.5\dots$ to $\nu_{135} = \nu_{270} = 1.124\dots \cdot 10^5$.

This is not a unique occurrence. It appears in many dimensions, most notably in those which are multiples of $24$. It is perhaps worth noting that the prime factorization of the number of minimal vectors in such a sequence (e.g. 222, 228, 234, 252, 270) all have the same number of distinct prime factors, and at each step at least one large prime factor is converted into lower prime factors. For instance, $222 = 2 \cdot 3 \cdot 37$ while $228 = 2^2 \cdot 3 \cdot 19$, which converts the $37$ to $2 \cdot 19$. This reduction of the largest term in the denominator of (\ref{nu_L}) drives up $\nu_n$ but holds $\omega(n)$ constant so drives down $\A_n$ as $n$ increases.

\begin{table} 
\centering
\begin{tabular}{|c|c|c|c|c|c|} \hline
$\phi(n)$ & $n$ & $\C_n$ & $\A_n$ & $\nu_n$ & $\Pi_n$ \\
\hline \hline \hline
6 & $7$ & {\bf 0.166...} & 0.166... & {\bf 1.666...} & 11.662... \\ \hline
6 & $9 = 3^2$ & 0.5 & {\bf 0.125} & 1.539...  & {\bf 18.475...} \\ \hline
\hline \hline
8 & $15 = 3 \cdot 5$ & 0.5 & 0.214... & {\bf 3.640...} & {\bf 31.857...} \\ \hline
8 & $16 = 2^4$ & {\bf 0} & {\bf 0} & 1 & -- \\ \hline
8 & $20 = 2^2 \cdot 5$ & 0.25 & 0.157... & 2.048  & 16.213... \\ \hline
8 & $24 = 2^3 \cdot 3$ & 0.5 & 0.090... & 1.777... & 29.333... \\ \hline
\hline \hline
24 & $35 = 5 \cdot 7$ & 0.25 & 0.088... & {\bf 66.194...} & 1094.055... \\ \hline
24 & $39 = 3 \cdot 13$ & 0.5 & 0.078... & 27.953... & 575.369... \\ \hline
24 & $45 = 3^2 \cdot 5$ & 0.5 & 0.068... & 48.263... & {\bf 1327.257...} \\ \hline
24 & $52 = 2^2 \cdot 13$ & {\bf 0.083...} & 0.04 & 4.975... & 134.741... \\ \hline
24 & $56 = 2^3 \cdot 7$ & 0.166... & 0.037... & 7.706... & 242.742...  \\ \hline
24 & $72 = 2^3 \cdot 3^2$ & 0.5 & {\bf 0.028...} & 5.618... & 294.979...  \\ \hline
24 & $84 = 2^2 \cdot 3 \cdot 7$ & 0.5 & 0.073... & 43.297... & 1035.542...  \\ \hline
\hline \hline
72 & $73$ & {\bf 0.013...} & 0.013... & 5.200... & 379.606... \\ \hline
72 & $91 = 7 \cdot 13$ & 0.166... & 0.033... & 56350.535... & $2.136... \cdot 10^6$  \\ \hline
72 & $95 = 5 \cdot 19$ & 0.25 & 0.031... & 32670.615... & $1.350... \cdot 10^6$  \\ \hline
72 & $111 = 3 \cdot 37$ & 0.5 & 0.027... & 2447.523... & $1.383... \cdot 10^5$  \\ \hline
72 & $117 = 3^2 \cdot 13$ & 0.5 & 0.025... & 21841.954... & $1.372... \cdot 10^6$  \\ \hline
72 & $135 = 3^3 \cdot 5$ & 0.5 & 0.022... & ${\bf 1.124... \cdot 10^5}$ & ${\bf 9.415... \cdot 10^6}$  \\ \hline
72 & $148 = 2^2 \cdot 37$ & 0.027... & 0.013... & 13.798... & 1035.267...  \\ \hline
72 & $152 = 2^3 \cdot 19$ & 0.055... & 0.013... & 51.545... & 4081.677...  \\ \hline
72 & $216 = 2^3 \cdot 3^3$ & 0.5 & {\bf 0.009...} & 177.376... & 28469.292...  \\ \hline
72 & $228 = 2^2 \cdot 3 \cdot 19$ & 0.5 & 0.026... & 9142.921... & $5.452... \cdot 10^5$ \\ \hline
72 & $252 = 2^2 \cdot 3^2 \cdot 7$ & 0.5 & 0.024. & 81171.032... & $5.918... \cdot 10^6$ \\ \hline
\hline \hline
160 & $187 = 11 \cdot 17$ & 0.1 & 0.016... & $1.163...  \cdot 10^9$ & $8.428... \cdot 10^{10}$ \\ \hline
160 & $205 = 5 \cdot 41$ & 0.25 & 0.014... & $3.928...  \cdot 10^8$ & $3.594... \cdot 10^{10}$\\ \hline
160 & $328 = 2^3 \cdot 41$ & {\bf 0.025} & 0.006... & $233.162...$ & 77912.090... \\ \hline
160 & $352 = 2^5 \cdot 11$ & 0.1 & 0.005... & $104646.972...$ & $2.014... \cdot 10^7$ \\ \hline
160 & $400 = 2^4 \cdot 5^2$ & 0.25 & {\bf 0.005...} & $1.684... \cdot 10^{6}$ & $4.191... \cdot 10^8$ \\ \hline
160 & $440 = 2^3 \cdot 5 \cdot 11$ & 0.25 & 0.013... & $1.763... \cdot 10^{11}$ & $1.769... \cdot 10^{13}$\\ \hline
160 & $492 = 2^2 \cdot 3 \cdot 41$ & 0.5 & 0.012... & $2.318... \cdot 10^{7}$ & $2.911... \cdot 10^9$ \\ \hline
160 & $528 = 2^4 \cdot 3 \cdot 11$ & 0.5 & 0.011... & $1.040... \cdot 10^{10}$ & $1.505... \cdot 10^{12}$\\ \hline
160 & $600 = 2^3 \cdot 3 \cdot 5^2$ & 0.5 & 0.010... & $1.675... \cdot 10^{11}$ & $3.131... \cdot 10^{13}$ \\ \hline
160 & $660 = 2^2 \cdot 3 \cdot 5 \cdot 11$ & 0.5 & 0.021... & ${\bf 1.753... \cdot 10^{16}}$ & $\bf{1.699... \cdot 10^{18}}$ \\ \hline
\end{tabular}
\medskip
\caption{Examples of coherence, average coherence, orthogonality defect and product measure values for cyclotomic lattices} 
\label{coh_orth_ex}
\end{table}

For comparison purposes, we also record here the values of coherence, average coherence, orthogonality defect and product measure for the standard irreducible root lattices. We start by briefly recalling some standard notation. A lattice is called \textit{irreducible} if it is not a direct sum of nonzero sublattices. A \textit{root} in a lattice is a vector of squared-norm equal to~$2$, and an irreducible lattice is called a \textit{root lattice} if it is generated by its roots. In this case, the roots are the minimal vectors of the lattice. There are precisely two infinite families of irreducible root lattices, denoted $A_n$ and $D_n$, as well as the three exceptional examples $E_6$, $E_7$ and $E_8$. We already defined $A_n$ in~\eqref{A_n}, and now recall that
\begin{equation}\label{D_n}
D_n = \left\{\bx \in \zed^n : \sum_{i = 1}^n x_i \in 2 \zed \right\},\ E_8 = D_8 \cup \left\{\frac{1}{2} \left(\sum_{i = 1}^8 \be_i\right) + D_8 \right\},
\end{equation}
where $\be_i$ are the standard basis vectors in the corresponding $\zed^n$. Additionally,
\begin{equation}\label{E_76}
E_7 = \left\{\bx \in E_8 : \left< \bx, \be_7+\be_8 \right> = 0 \right\},\ E_6 = \left\{\bx \in E_7 : \left< \bx, \be_6+\be_8 \right> = 0 \right\}.
\end{equation}
We refer the reader to~\cite{martinet} (Chapter~4) or~\cite{conway:sloane} (Chapter~4) for the detailed information on the properties of root lattices. We will mention that, due to the remarkable symmetry properties of root lattices, their minimal vectors are indistinguishable in the following sense. Let $L$ be a root lattice. Then for each vector $\bx \in S(L)$ there is the same number of vectors $\bwy \in S(L)$ that have nonzero inner product $\left< \bx, \bwy \right>$ (\cite{martinet}, Proposition 4.10.12). Standard integrality conditions limit the only other possible inner product value to $|\left< \bx, \bwy \right>| = 1$. With this in mind, the calculation of the average coherence of root lattices becomes straightforward, using Proposition~4.2.2 and Theorems~4.3.3, 4.4.4, 4.5.2 and~4.5.3 of~\cite{martinet}. The values held by the coherence, average coherence, orthogonality defect and product measure on the corresponding root lattices $A_n$ for $n \geq 2$, $D_n$ for $n \geq 4$, $E_6$, $E_7$, and $E_8$ are given in Table~\ref{root_table}.

\begin{table}[h]
\centering
\begin{tabular}{|c|c|c|c|c|c|} \hline
Lattice $L$ & $|S'(L)|$ & $\C(L)$ & $\A(L)$ & $\nu(L)$ & $\Pi(L)$ \\
\hline \hline \hline
$A_n$ & $\frac{n(n+1)}{2}$ & 0.5 & $\frac{2}{n+2}$ & $\frac{2^{\frac{n}{2}}}{n+1}$ & $(n+2) 2^{\frac{n-4}{2}}$ \\
 \hline
$D_n$ & $n(n-1)$ & 0.5 & $\frac{2(n-2)}{n^2 -n-1}$ & $2^{\frac{n - 4}{2}}$ & $\frac{(n-1)(n^2 -n-1)}{n-2} 2^{\frac{n-6}{2}}$ \\
 \hline
$E_6$ & 36 & 0.5 & $\frac{2}{7}$ & $\frac{8}{3}$ & 56 \\
 \hline
$E_7$ & 63 & 0.5 & $\frac{8}{31}$ & $4 \sqrt{2}$ & 13.138\dots \\
 \hline
$E_8$ & 120 & 0.5 & $\frac{28}{119}$ & 16 & 1020 \\
 \hline
\end{tabular}
\medskip
\caption{Coherence, average coherence, orthogonality defect and product measure values for root lattices}
\label{root_table}
\end{table}

This data suggests that root lattices are generally better than cyclotomic lattices at simultaneously minimizing average coherence and maximizing orthogonality defect, however are worse at minimizing maximal coherence. Indeed, suppose some large $p$ is the smallest prime dividing $n$ and let $d = \phi(n)$, then $\Lambda_{\que(\zeta_n)})$ is a lattice in $\real^d$ with maximal coherence $1/(p-1)$, while $A_d$ and $D_d$ are root lattices in the same dimension with maximal coherence~$1/2$.

In fact, an interesting feature of the cyclotomic lattices, in contrast with the root lattices, is that their maximal and average coherence are about the same on the average as $n \to \infty$. Indeed, $\C_n = 1/(\eta(n) -1)$, where $\eta(n)$ is the smallest prime divisor of $n$. Now, the average order of $\eta(n)$ is known to be $(1+o(1))n/ 2 \log n$ as $n \to \infty$ (see~\cite{kalecki}). Hence the average order of $\C_n$ is $\frac{2 \log n}{n}$. On the other hand, the average order of $\omega(n)$ is $\log \log n$ (see Theorem~430 of~\cite{hardy_wright}). Combining this observation with~\eqref{av_coh_formula}, we see that the average order of $\A_n$ is $\frac{\log 2 \log n}{n}$.

\subsection*{Acknowledgment} We thank the anonymous referee for many helpful remarks and suggestions that improved the quality of the paper. Fukshansky was partially supported by the Simons Foundation grant \#519058.

\EditInfo{December 14, 2020}{April 15, 2021}{Karl Dilcher}

\end{paper}